 
\magnification=\magstep1       
\hsize=5.9truein                     
\vsize=8.5truein                       
\parindent 0pt
\parskip=\smallskipamount
\mathsurround=1pt
\hoffset=.25truein                     
\voffset=2\baselineskip               
%
%
\def\today{\ifcase\month\or
  January\or February\or March\or April\or May\or June\or
  July\or August\or September\or October\or November\or December\fi
  \space\number\day, \number\year}
%
 at 10truept
%
\newcount\dispno      
\dispno=1\relax       
\newcount\refno       
\refno=1\relax        
\newcount\citations   
\citations=0\relax    
\newcount\sectno      
\sectno=0\relax       
\newbox\boxscratch    
%

%
%
%
\def\Section#1#2{\global\advance\sectno by 1\relax%
\label{Section\noexpand~\the\sectno}{#2}%
\smallskip
\goodbreak
\setbox\boxscratch=\hbox{\bf Section \the\sectno.~}%
{\hangindent=\wd\boxscratch\hangafter=1
\noindent{\bf Section \the\sectno.~#1}\nobreak\smallskip\nobreak}}
%
\def\sqr#1#2{{\vcenter{\vbox{\hrule height.#2pt
              \hbox{\vrule width.#2pt height#1pt \kern#1pt
              \vrule width.#2pt}
              \hrule height.#2pt}}}}
\def\square{$\mathchoice\sqr34\sqr34\sqr{2.1}3\sqr{1.5}3$}
\def\endproof{~~\hfill\square\par\medbreak}
\def\noproof{~~\hfill\square}
%
%
\def\proc#1#2#3{{\hbox{${#3 \subseteq} \kern -#1cm _{#2 /}\hskip 0.05cm $}}}
\def\propcont{\mathchoice\proc{0.17}{\scriptscriptstyle}{}
                         \proc{0.17}{\scriptscriptstyle}{}
                         \proc{0.15}{\scriptscriptstyle}{\scriptstyle }
                         \proc{0.13}{\scriptscriptstyle}{\scriptscriptstyle}}
%

%
\def\normalin{\hbox{\raise0.045cm \hbox
                   {$\underline{\triangleleft }$}\hskip0.02cm}}
%
%
\def\'#1{\ifx#1i{\accent"13 \i}\else{\accent"13 #1}\fi}
%
%
%
\def\semidirect{\rlap{$\times$}\kern+7.2778pt \vrule height4.96333pt
width.5pt depth0pt\relax\;}
%
%
\def\prop#1#2{\noindent{\bf Proposition~\the\sectno.\the\dispno. }%
\label{Proposition\noexpand~\the\sectno.\the\dispno}{#1}\global\advance\dispno 
by 1{\it #2}\smallbreak}
\def\thm#1#2{\noindent{\bf Theorem~\the\sectno.\the\dispno. }%
\label{Theorem\noexpand~\the\sectno.\the\dispno}{#1}\global\advance\dispno
by 1{\it #2}\smallbreak}
\def\cor#1#2{\noindent{\bf Corollary~\the\sectno.\the\dispno. }%
\label{Corollary\noexpand~\the\sectno.\the\dispno}{#1}\global\advance\dispno by
1{\it #2}\smallbreak}
\def\defn{\noindent{\bf
Definition~\the\sectno.\the\dispno. }\global\advance\dispno by 1\relax}
\def\lemma#1#2{\noindent{\bf Lemma~\the\sectno.\the\dispno. }%
\label{Lemma\noexpand~\the\sectno.\the\dispno}{#1}\global\advance\dispno by
1{\it #2}\smallbreak}
\def\rmrk#1{\noindent{\bf Remark~\the\sectno.\the\dispno.}%
\label{Remark\noexpand~\the\sectno.\the\dispno}{#1}\global\advance\dispno
by 1\relax}
\def\proof{\noindent{\it Proof: }}
\def\numbeq#1{\the\sectno.\the\dispno\label{\the\sectno.\the\dispno}{#1}%
\global\advance\dispno by 1\relax}

\def\comm#1,#2{\left[#1{,}#2\right]}
\newdimen\boxitsep \boxitsep=0 true pt
\newdimen\boxith \boxith=.4 true pt 
\newdimen\boxitv \boxitv=.4 true pt
\gdef\boxit#1{\vbox{\hrule height\boxith
                    \hbox{\vrule width\boxitv\kern\boxitsep
                          \vbox{\kern\boxitsep#1\kern\boxitsep}%
                          \kern\boxitsep\vrule width\boxitv}
                    \hrule height\boxith}}
\def\square{\ \hbox{\vrule height7.5pt depth1.5pt width 6pt}\par}
\outer\def\square{\ifmmode\else\hfill\fi
   \setbox0=\hbox{} \wd0=6pt \ht0=7.5pt \dp0=1.5pt
   \raise-1.5pt\hbox{\boxit{\box0}\par}
}

\def\frac#1/#2{\leavevmode\kern.1em
              \raise.5ex\hbox{\the\scriptfont0 #1}\kern-.1em
              /\kern\.15em\lower.25ex\hbox{\the\scriptfont0 #2}}
\def\incnoteq{\lower.1ex \hbox{\rlap{\raise 1ex
     \hbox{$\scriptscriptstyle\subset$}}{$\scriptscriptstyle\not=$}}}
%
%


\def\propcontup{\bigcup\!\!\!\rlap{\kern+.2pt$\backslash$}\,\kern+1pt\vert}
%
%
%
\def\label#1#2{\immediate\write\aux%
{\noexpand\def\expandafter\noexpand\csname#2\endcsname{#1}}}
%
\def\ifundefined#1{\expandafter\ifx\csname#1\endcsname\relax}
%
%
\def\ref#1{%
\ifundefined{#1}\message{! No ref. to #1;}%
 \else\csname #1\endcsname\fi}
%
%
\def\refer#1{%
\the\refno\label{\the\refno}{#1}%
\global\advance\refno by 1\relax}
%
%
\def\cite#1{%
\expandafter\gdef\csname x#1\endcsname{1}%
\global\advance\citations by 1\relax
\ifundefined{#1}\message{! No ref. to #1;}%
\else\csname #1\endcsname\fi}
%
%
\font\bb=msbm10 
 at 8truept      
%
%
%

\def\Z{\hbox{\bb Z}}

\def\Z{\hbox{\bb Z}}                     

\newread\aux
\immediate\openin\aux=\jobname.aux
\ifeof\aux \message{! No file \jobname.aux;}
\else \input \jobname.aux \immediate\closein\aux \fi
\newwrite\aux
\immediate\openout\aux=\jobname.aux

\font\smallheadfont=cmr8 at 8truept

\headline={\ifnum\pageno<2{\hfill}\else{\ifodd\pageno\rightheadline
\else\leftheadline\fi}\fi}
\def\leftheadline{\smallheadfont A. Magidin\hfil}
\def\rightheadline{\hfil\smallheadfont A generalized argument for
dominions in varieties of groups\quad\quad}
 
\centerline{\bf A generalized argument for dominions in varieties of groups}
\centerline{Arturo Magidin\footnote*{The author was
supported in part by a fellowship from the Programa de Formaci\'on y
Superaci\'on del Personal Acad\'emico de la UNAM, administered by the
DGAPA.}}
\smallskip
{\parindent=20pt
\narrower\narrower
\noindent{\smallheadfont{Abstract. An argument used to
show that certain varieties of nilpotent groups have instances of
nontrivial dominions is considered, and generalized.
The same is done with the argument used to show that there are nontrivial
dominions in the variety of metabelian groups, to suggest how this
general technique may be used.\par}}}
\bigskip
\medskip
 
\footnote{}{\noindent\smallheadfont Mathematics Subject
Classification:
08B25, 20E10 (primary)}
\footnote{}{\noindent\smallheadfont Keywords:dominion}

\Section{Introduction}{intro}
 
Let~$\cal C$ be a full subcategory of the category of all algebras (in
the sense of Universal Algebra) of a fixed type, which is closed under
passing to subalgebras.  Let $A\in {\cal C}$, and let~$B$ be a
subalgebra of~$A$. Recall that, in this situation, Isbell (see~{\bf
[\cite{isbellone}]}) defines the {\it dominion of~$B$ in~$A$} (in the
category ${\cal C}$) to be the intersection of all equalizer
subalgebras of~$A$ containing~$B$. Explicitly,
$${\rm dom}_A^{\cal C}(B)=\Bigl\{a\in A\bigm| \forall C\in {\cal C},\;
\forall f,g\colon A\to C,\ {\rm if}\ f|_B=g|_B{\rm\ then\ }
f(a)=g(a)\Bigr\}.$$

Note that ${\rm dom}_A^{\cal C}(B)$ always contains~$B$. If $B$ is
properly contained in its dominion, we will say that the dominion
of~$B$ in~$A$ is {\it nontrivial}, and call it {\it trivial\/} if it
equals~$B$. A category~${\cal C}$ {\it has instances nontrivial
dominions} if there is an algebra $A\in {\cal C}$, and a subalgebra
$B$ of~$A$ such that the dominion of~$B$ in~$A$ (in the
category~${\cal C}$) is~nontrivial.

In this work we will restrict our attention to the case where ${\cal
C}$ is a variety of groups. For the basic properties of dominions in
the context of varieties of groups, we direct the reader to~{\bf
[\cite{domsmetabprelim}]}. We recall the most important properties:
for a group~$G$, ${\rm dom}_G^{\cal C}(-)$ is a closure operator on
the lattice of subgroups of~$G$, and normal subgroups are
dominion-closed; dominions respect finite direct powers; and dominions
respect quotients. That is, if $G\in {\cal C}$, $H$ is a subgroup
of~$G$, and $N$ is a normal subgroup of~$G$ contained in~$H$, then
$${\rm dom}_{(G/N)}^{\cal C}\bigl(H/N\bigr) = \Bigl({\rm dom}_G^{\cal
C}(H)\Bigr)\bigm/ N.$$

All groups will be written multiplicatively unless otherwise
noted. All maps will be assumed to be group morphisms unless
otherwise noted. Given a group~$G$, we will denote the identity
element of~$G$ by~$e_G$, although we will omit the subscript if it is
clear from context. Given a group~$G$, and elements $x$ and~$y$
of~$G$, we denote their commutator by $[x,y]=x^{-1}y^{-1}xy$.

Given two
subsets $A,B$ of~$G$ (not necessarily subgroups), we denote by $[A,B]$
the subgroup of~$G$ generated by all elements $[a,b]$ with $a\in A$
and $b\in B$.  We also define inductively the left-normed commutators
of weight~\hbox{$c+1$}:
$$[x_1,\ldots,x_c,x_{c+1}] =
\bigl[[x_1,\ldots,x_c],x_{c+1}\bigr];\quad c\geq 2.$$
We will denote the center of~$G$ by $Z(G)$.
 
The following lemma, which is easily established by direct
computation will be useful in
subsequent considerations.
 
\lemma{commident}{The
following hold for any elements $x$,~$y$, $z$, and~$w$  of
an arbitrary group~$G$:}
{\parindent=30pt
\it
\item{(a)}{$xy = yx[x,y];\qquad x^y=x[x,y]$.}
\item{(b)}{$[x,y]^{-1} = [y,x]$.}
\item{(c)}{$[xy,z] = [x,z]^y[y,z] = [x,z][x,z,y][y,z]$.}
\item{(d)}{$[x,zw] = [x,w][x,z]^w = [x,w][x,z][x,z,w]$.\noproof}\par}

Varieties will be donoted by caligraphic letters, ${\cal A}$, ${\cal
V}$, ${\cal W}$, etc. We will denote the variety of all groups
by~${\cal G}$, and the variety consisting only of the trivial group
by~${\cal E}$. We will denote the variety of nilpotent groups of
class~two by~${\cal N}_2$, and the variety of metabelian groups (that
is, groups which are an extension of an abelian group by an abelian
group) by~${\cal A}^2$.

In \ref{generalargument} we will consider the argument used in~{\bf
[\cite{nildomsprelim}]} to prove the existence of nontrivial dominions
in~${\cal N}_2$ and some of its subvarieties, and we will modify it to
a more general context. The ideas can easily be modified to deal with
other arguments, and we do this in \ref{othergens}, where we
consider the argument used in~{\bf [\cite{domsmetabprelim}]} to
establish the existence of nontrivial dominions in~${\cal A}^2$.

The contents of this work are part of the author's doctoral
dissertation, which was conducted under the direction of Prof.~George
M.~Bergman, at the University of California at~Berkeley. It is my very
great pleasure to express my deep gratitude and indebtedness to
Prof.~Bergman for his advice and encouragement; the results in this
paper are the direct consequence of his suggestions, which have also
improved the work in ways too numerous to list explicitly; he also
helped me correct many mistakes. Any errors that remain, however, are
my own~responsibility.

\Section{Generalizing the $[x,y]^p$ argument}{generalargument}

In {\bf [\cite{nildomsprelim}]}, we studied subvarieties of~${\cal
N}_2$, and exhibited a large family of such subvarieties that had
nontrivial dominions. The basic idea was as follows: Given a group~$G$
in such a variety, and elements $x$ and~$y$ in~$G$, we looked
at~$[x,y]$. Since $G$ is nilpotent of class~two, the element $[x,y]$
commutes with both $x$ and~$y$. It follows that the identity
$[x^n,y]=[x,y]^n$ holds in~$G$ for every $n\in \Z$. If the commutator
subgroup of~$G$ does not have squarefree exponent, then we choose a
prime~$p$ whose squre divides the order of~$[x,y]$. Finally, we
let~$H$ be generated by $x^p$ and $y^p$. Then the elements $[x,y]^p$
lies in the dominion of~$H$. Some mild conditions on~$G$ guarantee
that $[x,y]^p$ does not lie in~$H$, and this gives an instance of a
nontrivial dominion; see~{\bf
[\cite{nildomsprelim}]}, Section~6for the~details..

In this section, we will generalize this process by replacing $x$ and
$y$ with words $v$ and~$w$, and replacing ${\cal N}_2$ with a variety
in which the commutator of $v$ and~$w$ commutes with both~$v$ and~$w$.

Let $\cal V$ be a variety of groups, let $v$ and~$w$ be two words
in~$m$ letters, and let $F$ be the relatively free ${\cal V}$-group in
$m$ generators, generated by $z_1,\ldots,z_m$.

\lemma{ifFthenalways}{Assume the notation of the preceding paragraph, and let
$${\bf z}=(z_1,\ldots,z_m).$$
Suppose~that in~$F$ the following identities hold:
$$\eqalign{[[v({\bf z}),w({\bf z})],v({\bf z})] &= e\cr
[[v({\bf z}),w({\bf z})],w({\bf z})] &= e.\cr}\leqno(\numbeq{displayone})$$
Then for any group~$G\in {\cal V}$, and any elements
$g_1,\ldots,g_m\in G$, the subgroup generated by $v(g_1,\ldots,g_m)$
and $w(g_1,\ldots,g_m)$ is nilpotent of class at most~two.}

\proof This follows because 
$$\langle x,y\,|\,[x,y,x]=[x,y,y]=e\rangle$$ is a presentation for the
relatively free ${\cal N}_2$-group of rank~two,
and hence the subgroup of~$F$ generated by
$v({\bf z})$ and~$w({\bf z})$ is a quotient of an ${\cal N}_2$-group.
Since the relations in question hold in the relatively free group of
the variety, it follows that they are inherited by all groups in the
variety. \endproof
 
For the remainder of this section, we will assume
that~(\ref{displayone}) holds in~${\cal V}$. Note
that~(\ref{displayone}) implies that the commutator subgroup
of~$\langle v({\bf z}),w({\bf z})\rangle$ is cyclic and generated
by~$[v({\bf z}),w({\bf z})]$.
 
Assume that $[v({\bf z}),w({\bf z})]$ is not of exponent~$k$ for
any square free positive integer $k>1$. We note in particular that
$[v({\bf z}),w({\bf z})]$ cannot be trivial. Let $k_0$ be the order of
$[v({\bf z}),w({\bf z})]$ in~$F$, if this order is
finite, and $k_0=0$ if $[v({\bf z}),w({\bf z})]$ has
infinite order. By assumption, $k_0$ is not square free. 
 
\def\vv(#1){v({\bf #1})}
\def\ww(#1){w({\bf #1})}
\def\zz(#1){z({\bf #1})}

Let $a_0$ be the order of $\vv(z)$ if the
latter is finite, and $a_0=0$ if the order is infinite, and
$b_0$ the order of $\ww(z)$ if this is finite, and $b_0=0$ if the order
of $\ww(z)$ is infinite. Note that we necessarily have $k_0|{\rm
gcd}(a_0,b_0)$.

Therefore, the group $G=\langle \vv(z),\ww(z)\rangle$ is a quotient of
the group $$K=\langle\, x,y\,\mid\,x^{a_0} = y^{b_0} =
[x,y]^{k_0} = [[x,y],x] = [[x,y],y] =
e\,\rangle.\leqno(\numbeq{defofK})$$
 
Thus, every element $g\in G$ can be written in the form
$$g=\vv(z)^a\ww(z)^b[\vv(z),\ww(z)]^c\leqno(\numbeq{shape})$$
where $0\leq a<a_0$ (if $a_0>0$, and $a\in\Z$ if $a_0=0$), $0\leq
b<b_0$ (if $b_0>0$, and $b\in\Z$ if $b_0=0$), and $0\leq c<k_0$ (if
$k_0>0$, and $c\in\Z$ if $k_0=0$). 

Suppose further that ${\bf z}$ is a
disjoint union of two sets of indeterminates
$$\eqalign{{\bf x}&=(x_1,\ldots,x_{n_1})\cr
           {\bf y}&=(y_1,\ldots,y_{n_2})\cr}$$
and that the word $v$ involves
only one of these subsets and $w$ involves only the other, say ${\bf
x}$ and ${\bf y}$ respectively. Under these extra assumptions, we
claim that the expression in (\ref{shape}) is unique.  Indeed, suppose
that we had a relation
$$\vv(x)^a\ww(y)^b[\vv(x),\ww(y)]^c=e.\leqno(\numbeq{relation})$$
We want to show that $a_0|a$, $b_0|b$ and $k_0|c$.
 
Consider the endomorphism $\psi$ of $F$ given by sending
$y_1,\ldots,y_{n_2}$ to
$e$ and leaving the $x_i$ unchanged. Applying $\psi$ to
(\ref{relation}), we obtain that $\vv(x)^a=e$. Therefore,
$a_0|a$. Applying the endomorphism of~$F$ that sends $x_1,\ldots,x_{n_1}$
to $e$ and leaves the $y_j$ unchanged we obtain, again from
(\ref{relation}), that $\ww(y)^b=e$, and hence $b_0|b$. Combining these
two relations with (\ref{relation}) we get that $[\vv(x),\ww(y)]^c=e$,
hence $k_0|c$. This proves our claim.

The above situation allows us to deduce that certain elements must lie
in the dominion of a given subgroup. Namely:

\lemma{halfnontrivial}{Let ${\cal V}$ be a variety of groups, and let
$\vv(x)$ and $\ww(y)$ be two words. Suppose that in~${\cal V}$, the
relatively free group~$F$ on $n_1+n_2$ generators $x_1,\ldots,x_{n_1}$
and $y_1,\ldots,y_{n_2}$ satisfies
the two identities~(\ref{displayone}). If $G\in {\cal V}$,
$g_1,\ldots,g_{n_1+n_2}\in G$, and~$H$ is
a subgroup of~$G$ which contains
$$v(g_1,\ldots,g_{n_1})^n,\qquad w(g_{n_1+1},\ldots,g_{n_1+n_2})^n$$
for some $n\in\Z$, then ${\rm dom}_G^{\cal V}(H)$ also contains
$$[v(g_1,\ldots,g_{n_1}),w(g_{n_1+1},\ldots,g_{n_1+n_2})]^n.$$}

\proof From \ref{ifFthenalways} it follows that all groups in~${\cal
V}$ satisfy the identities 
$$[a^n,b] = [a,b]^n = [a,b^n]$$
for all elements $a$ and~$b$, and for all $n\in Z$.

Write ${\bf g_1}=(g_1,\ldots,g_{n_1})$, and ${\bf
g_2}=(g_{n_1+1},\ldots,g_{n_1+n_2})$. 

Let $K\in {\cal V}$, and let $f,h\colon G\to K$ be two morphisms such
that $f|_H=h|_H$. Then
$$\eqalign{f\left([\vv(g_1),\ww(g_2)]^n\right) &=
                                f\left([\vv(g_1)^n,\ww(g_2)]\right)\cr 
                                &=
              \Bigl[f\bigl(\vv(g_1)^n\bigr),f\bigl(\ww(g_2)\bigr)\Bigr]\cr
    &= \Bigl[h\bigl(\vv(g_1)^n\bigr),f\bigl(\ww(g_2)\bigr)\Bigr]\cr
    &= \Bigl[h\bigl(\vv(g_1)\bigr)^n,f\bigl(\ww(g_2)\bigr)\Bigr]\cr
    &= \Bigl[h\bigl(\vv(g_1)\bigr),f\bigl(\ww(g_2)\bigr)\Bigr]^n\cr}$$
and by a symmetric argument, this term equals $h([\vv(g_1),\ww(g_2)]^n)$. In
particular, $[\vv(g_1),\ww(g_2)]^n$ lies in ${\rm dom}_{F}^{{\cal
V}}(H)$.\endproof

Note that in \ref{halfnontrivial} we have made no claims on whether
$[\vv(g_1),\ww(g_2)]^n$ lies in~$H$ or~not.

Still under the assumptions of the two paragraphs preceding
\ref{halfnontrivial}, let $p$ be
a prime such that $p^2|k_0$, and let~$H$ be
the subgroup of~$G$ generated by $\vv(x)^p$ and $\ww(y)^p$. The
elements of~$H$ correspond to those elements of~$G$ that can be written as
in~(\ref{shape}), with $p|a$, $p|b$, and~$p^2|c$.
 
By \ref{halfnontrivial}, we know that
$[\vv(x),\ww(y)]^p$ lies in
${\rm dom}_{F}^{{\cal V}}({H})$. Since elements of~$G$ can be written
uniquely in the form given in~(\ref{shape}), $[\vv(x),\ww(y)]^p$ does
not lie in the subgroup~$H$. This gives the following result:
 
\thm{firstnontrivial}{Let $\cal V$ be a variety of groups, and
$v(x_1,\ldots,x_{n_1})$, $w(y_1,\ldots,y_{n_2})$ be two words. 
Suppose that in $\cal V$, the relatively free group $F$ on $n_1+n_2$ 
generators, $x_1,\ldots,x_{n_1}$, $y_1,\ldots,y_{n_2}$ satisfies
$$\eqalign{[[\vv(x),\ww(y)],\vv(x)]&=e\cr
           [[\vv(x),\ww(y)],\ww(y)]&=e.\cr}$$
Let $k$ denote the order of $[\vv(x),\ww(y)]$ if this is finite, and
$k=0$ if the order is infinite. Suppose $k$ is divisible by the square
of some prime number~$p$.
Then there are instances of nontrivial dominions in
$\cal V$. Namely, for $H = \langle \vv(x)^p,\ww(y)^p\rangle$,
we have $H\propcont {\rm dom}_{F}^{{\cal V}}({H})$.\noproof}

Recall that if ${\cal V}$ and ${\cal W}$ are varieties of groups,
corresponding to the fully invariant subgroups $V$ and~$W$
of~$F_{\infty}$, respectively, we define the
variety $[{\cal V},{\cal W}]$ to be the variety corresponding to the
fully invariant subgroup $[V,W]$ of~$F_{\infty}$. It is not hard to
verify that $[{\cal V},{\cal W}]$ is defined by the identities
$[v,w]$, where $v$ is an identity of~${\cal V}$, and $w$ an identity
of~${\cal W}$.

Also, if we let ${\cal E}$ denote the variety consisting only of the
trivial group, then for a given variety ${\cal V}$ we define the
variety of {\it center-by-${\cal V}$ groups} to be the variety $[{\cal
E},{\cal V}]$. These groups can be described as groups 
$G$ such that $G/Z(G)\in {\cal V}$, where $Z(G)$ denotes the
center of~$G$.

As such, it is not hard to verify that the
center-by-abelian groups are the nilpotent groups of class~$2$ (that
is, $[{\cal E},{\cal A}]={\cal N}_2$), and in general that $[{\cal
E},{\cal N}_c]={\cal N}_{c+1}$. Thus,
\ref{firstnontrivial} (together with the well-known fact that the
relatively free groups in~${\cal N}_c$ are torsion-free) could be used
to prove a result from~{\bf [\cite{nildomsprelim}]},
that there are nontrivial dominions in the variety of
nilpotent groups of class~$c>1$.

Also, it is not hard to verify that given a variety ${\cal V}$,
$[{\cal V},{\cal V}]$ is the variety of abelian-by-${\cal V}$ groups,
that is the varieyt ${\cal AV}$; and that the variety $[[{\cal
V},{\cal V}],{\cal V}]$ is the variety of ${\cal N}_2$-by-${\cal V}$
groups, ${\cal N}_2{\cal V}$. This variety satisfies the hypothesis
of~\ref{firstnontrivial}, so we obtain the following:

\cor{ntwobyv}{If ${\cal V}$ is any variety of groups, and ${\cal
V}\not={\cal G}$, then the variety ${\cal N}_2{\cal V}$ has instances of
nontrivial dominions.\noproof}

\Section{Other generalizations}{othergens}

The ideas of \ref{generalargument} can be expanded to generalize other
arguments used to prove nontriviality of dominions in certain varieties.
In this section, we will
generalize the argument in~{\bf [\cite{domsmetabprelim}]} as an indication of
how this would be~done.

In {\bf [\cite{domsmetabprelim}]} we proved that given a group $G\in
{\cal A}^2$, the variety of metabelian groups, and a subgroup~$H$
of~$G$, if $x\in H\cap [G,G]$, and $y,z\in G$ are such that $[x,y]$
and $[x,z]$ lie in~$H$, then $[x,y,z]$ lies in~${\rm dom}_G^{{\cal
A}^2}(H)$. See Lemma~3.25 in~{\bf [\cite{domsmetabprelim}]} for
the details. We now proceed as in the previous section to lay the
groundwork for the~generalization.

\lemma{berggenforelements}{Let $G$ be a group, and $u,v,w\in G$. If
$$\bigl[v,[w^{-1},z^{-1}]\bigr] = \bigl[[v,w],[v,z]\bigr]=e$$
then $[v,w,z]=[v,z,w]$.}

\proof Since $v$ commutes with $[w^{-1},z^{-1}]$, we have that
$$v^{wzw^{-1}z^{-1}}= v$$
and therefore,
$v^{wz} = v^{zw}$.
Left-multiplying by $v^{-1}$ we get
$[v,wz] = [v,zw]$.
Expanding each side using commutator identities, we have
$$[v,w][v,z][v,w,z] = [v,z][v,w][v,z,w]$$
and since $[v,w]$ commutes with $[v,z]$ by hypothesis, we may cancel
the first two terms from both sides to get
$$[v,w,z] = [v,z,w]$$
as~claimed.\endproof

\lemma{firststepberggen}{Suppose that a variety of groups~${\cal V}$
satisfies the identities
$$\eqalign{\bigl[v({\bf x}), [w({\bf y})^{-1}, z({\bf x})^{-1}]\bigr]&= e\cr
\bigl[ [v({\bf x}),w({\bf y})], [v({\bf x}),z({\bf x})]\bigr]
&=e\cr}\leqno(\numbeq{generalizingmetab})$$ for certain
group-theoretic words $v$, $w$ and~$z$ in tuples of variables~${\bf
x}$ and~${\bf y}$. Then ${\cal V}$ also satisfies the identity
$$[v({\bf x}),w({\bf y}),z({\bf x})] = [v({\bf x}), z({\bf x}), w({\bf
y})].$$}

\proof This follows from \ref{berggenforelements}.\endproof

Now we may apply the argument from~{\bf [\cite{domsmetabprelim}]} to~prove:

\thm{generalizedomsinmetab}{Let ${\cal V}$ be a variety and suppose
that ${\cal V}$ satisfies the identities~(\ref{generalizingmetab}). If
$G\in {\cal V}$, then for any particular tuple~${\bf x}$ of elements
of~$G$, the dominion of the subgroup generated by $v({\bf x})$,
$[v({\bf x}),w({\bf x})]$ and~$[v({\bf x}),z({\bf x})]$ contains
$$[v({\bf x}),w({\bf x}),z({\bf x})].$$}

\proof Let $H$ be the subgroup of~$G$ generated by the elements $\vv(x)$,
$[\vv(x),\ww(x)]$, and~$[\vv(x),\zz(x)]$. Let $K\in {\cal V}$, and let
$f,g\colon G\to K$ be two group morphisms such that $f|_H=g|_H$. Then:
$$\eqalignno{f\Bigl([\vv(x),\ww(x),\zz(x)]\Bigr)&=
\Bigl[f\bigl([\vv(x),\ww(x)]\bigr),f(\zz(x))\Bigr]\cr
&=\Bigl[g\bigl([\vv(x),\ww(x)]\bigr),f(\zz(x))\Bigr]\qquad
\hbox{(since $[\vv(x),\ww(x)]\in H$)}\cr
&=\Bigl[g(\vv(x)),g(\ww(x)),f(\zz(x))\Bigr]\cr
&=\Bigl[f(\vv(x)),f(\zz(x)),g(\ww(x))\Bigr]\cr
&\qquad\qquad \hbox{(using
\ref{firststepberggen} and the fact that $\vv(x)\in H$)}\cr
&=\Bigl[f\bigl([\vv(x),\zz(x)]\bigr),g(\ww(x))\Bigr]\cr
&=\Bigl[g\bigl([\vv(x),\zz(x)]\bigr),g(\zz(x))\Bigr]\qquad\quad \hbox{(since
$[\vv(x),\zz(x)]\in H$)}\cr
&=g\Bigl([\vv(x),\zz(x),\ww(x)]\Bigr)\cr
&=g\Bigl([\vv(x),\ww(x),\zz(x)]\Bigr)\qquad\qquad\quad \hbox{(by
\ref{firststepberggen})}\cr}$$
so $[\vv(x),\ww(x),\zz(x)]\in {\rm dom}_G^{\cal V}(H)$, as~claimed.\endproof

From this, we can deduce the following strengthening of Lemma~3.25
in~{\bf [\cite{domsmetabprelim}]}:

\cor{domsinad}{Let $d\geq 2$, and let ${\cal A}^d$ be the variety of
all solvable groups of solvability length at most~$d$. Let $G\in {\cal
A}^d$, and let $H$ be a subgroup of~$G$. If $x\in G^{(d-1)}$, $y,z\in
G^{(d-2)}$ are such that $x$, $[x,y]$, and $[x,z]$ lie in~$H$, then
$[x,y,z]\in {\rm dom}_G^{{\cal A}^d}(H)$.\noproof}

Similar arguments may be used to generalize other results about
specific varieties, which are done through word-theoretic arguments
such as the~above. A related general result on dominions in the
context of Universal Algebra can be found
in~{\bf [\cite{wordsprelim}]}.

%
\ifnum0<\citations{\par\bigbreak
\filbreak{\bf References}\par\frenchspacing}\fi
%
\ifundefined{xthreeNB}\else
\item{\bf [\refer{threeNB}]}{Baumslag, G{.}, Neumann, B{.}H{.},
Neumann, H{.}, and Neumann, P{.}M. {\it On varieties generated by a
finitely generated group.\/} {\sl Math.\ Z.} {\bf 86} (1964)
pp.~\hbox{93--122}. {MR:30\#138}}\par\filbreak\fi
\ifundefined{xbergman}\else
\item{\bf [\refer{bergman}]}{Bergman, George M. {\it An Invitation to
General Algebra and Universal Constructions.\/} {\sl Berkeley
Mathematics Lecture Notes 7\/} (1995).}\par\filbreak\fi
\ifundefined{xordersberg}\else
\item{\bf [\refer{ordersberg}]}{Bergman, George M. {\it Ordering
coproducts of groups and semigroups.\/} {\sl J. Algebra} {\bf 133} (1990)
no. 2, pp.~\hbox{313--339}. {MR:91j:06035}}\par\filbreak\fi
\ifundefined{xbirkhoff}\else
\item{\bf [\refer{birkhoff}]}{Birkhoff, Garrett. {\it On the structure
of abstract algebras.\/} {\sl Proc.\ Cambridge\ Philos.\ Soc.} {\bf
31} (1935), pp.~\hbox{433--454}.}\par\filbreak\fi
\ifundefined{xbrown}\else
\item{\bf [\refer{brown}]}{Brown, Kenneth S. {\it Cohomology of
Groups, 2nd Edition.\/} {\sl Graduate texts in mathematics 87\/},
Springer Verlag,~1994. {MR:96a:20072}}\par\filbreak\fi
\ifundefined{xmetab}\else
\item{\bf [\refer{metab}]}{Golovin, O. N. {\it Metabelian products of
groups.\/}
{\sl American Mathematical Society Translations}, series 2, {\bf 2} (1956),
pp.~\hbox{117--131.} {MR:17,824b}}\par\filbreak\fi
\ifundefined{xhall}\else
\item{\bf [\refer{hall}]}{Hall, M. {\it The Theory of Groups.\/}
Mac~Millan Company,~1959. {MR:21\#1996}}\par\filbreak\fi
\ifundefined{xphall}\else
\item{\bf [\refer{phall}]}{Hall, P. {\it Verbal and marginal
subgroups.} {\sl J.\ Reine\ Angew.\ Math.\/} {\bf 182} (1940)
pp.~\hbox{156--157.} {MR:2,125i}}\par\filbreak\fi
\ifundefined{xheineken}\else
\item{\bf [\refer{heineken}]}{Heineken, H. {\it Engelsche Elemente der
L\"ange drei,\/} {\sl Illinois Journal of Math.} {\bf 5} (1961)
pp.~\hbox{681--707.} {MR:24\#A1319}}\par\filbreak\fi
\ifundefined{xherman}\else
\item{\bf [\refer{herman}]}{Herman, Krzysztof. {\it Some remarks on
the twelfth problem of Hanna Neumann.\/} {\sl Publ.\ Math.\ Debrecen}
{\bf 37} (1990)  no. 1--2, pp.~\hbox{25--31.} {MR:91f:20030}}\par\filbreak\fi
\ifundefined{xherstein}\else
\item{\bf [\refer{herstein}]}{Herstein, I.~N. {\it Topics in
Algebra.\/} Blaisdell Publishing Co.,~1964.}\par\filbreak\fi
\ifundefined{xepisandamalgs}\else
\item{\bf [\refer{episandamalgs}]}{Higgins, Peter M. {\it Epimorphisms
and amalgams.} {\sl
Colloq.\ Math.} {\bf 56} no.~1 (1988) pp.~\hbox{1--17.}
{MR:89m:20083}}\par\filbreak\fi
\ifundefined{xhigmanpgroups}\else
\item{\bf [\refer{higmanpgroups}]}{Higman, Graham. {\it Amalgams of
$p$-groups.\/} {\sl J. of~Algebra} {\bf 1} (1964)
pp.~\hbox{301--305.} {MR:29\#4799}}\par\filbreak\fi
\ifundefined{xhigmanremarks}\else
\item{\bf [\refer{higmanremarks}]}{Higman, Graham. {\it Some remarks
on varieties of groups.\/} {\sl Quart.\ J.\ of Math.\ (Oxford) (2)} {\bf
10} (1959), pp.~\hbox{165--178.} {MR:22\#4756}}\par\filbreak\fi
\ifundefined{xhughes}\else
\item{\bf [\refer{hughes}]}{Hughes, N.J.S. {\it The use of bilinear
mappings in the classification of groups of class~$2$.\/} {\sl Proc.\
Amer.\ Math.\ Soc.\ } {\bf 2} (1951) pp.~\hbox{742--747.}
{MR:13,528e}}\par\filbreak\fi
\ifundefined{xisbelltwo}\else
\item{\bf [\refer{isbelltwo}]}{Howie, J.~M., Isbell, J.~R. {\it
Epimorphisms and dominions II.\/} {\sl Journal of Algebra {\bf
6}}(1967) pp.~\hbox{7--21.} {MR:35\#105b}}\par\filbreak\fi
\ifundefined{xisaacs}\else
\item{\bf [\refer{isaacs}]}{Isaacs, I.M., Navarro, Gabriel. {\it
Coprime actions, fixed-point subgroups and irreducible induced
characters.} {\sl J.~of Algebra} {\bf 185} (1996) no.~1,
pp.~\hbox{125--143.} {MR:97g:20009}}\par\filbreak\fi
\ifundefined{xisbellone}\else
\item{\bf [\refer{isbellone}]}{Isbell, J. R. {\it Epimorphisms and
dominions} in {\sl 
Proc.~of the Conference on Categorical Algebra, La Jolla 1965,\/}
pp.~\hbox{232--246.} Lange and Springer, New
York~1966. MR:35\#105a (The statement of the
Zigzag Lemma for {\it rings} in this paper is incorrect. The correct
version is stated in~{\bf [\cite{isbellfour}]})}\par\filbreak\fi
\ifundefined{xisbellthree}\else
\item{\bf [\refer{isbellthree}]}{Isbell, J. R. {\it Epimorphisms and
dominions III.} {\sl Amer.\ J.\ Math.\ }{\bf 90} (1968)
pp.~\hbox{1025--1030.} {MR:38\#5877}}\par\filbreak\fi
\ifundefined{xisbellfour}\else
\item{\bf [\refer{isbellfour}]}{Isbell, J. R. {\it Epimorphisms and
dominions IV.} {\sl Journal\ London Math.\ Society~(2),}
{\bf 1} (1969) pp.~\hbox{265--273.} {MR:41\#1774}}\par\filbreak\fi
\ifundefined{xjones}\else
\item{\bf [\refer{jones}]}{Jones, Gareth A. {\it Varieties and simple
groups.\/} {\sl J.\ Austral.\ Math.\ Soc.} {\bf 17} (1974)
pp.~\hbox{163--173.} {MR:49\#9081}}\par\filbreak\fi
\ifundefined{xjonsson}\else
\item{\bf [\refer{jonsson}]}{J\'onsson, B. {\it Varieties of groups of
nilpotency three.} {\sl Notices Amer.\ Math.\ Soc.} {\bf 13} (1966)
pp.~488.}\par\filbreak\fi
\ifundefined{xwreathext}\else
\item{\bf [\refer{wreathext}]}{Kaloujnine, L. and Krasner, Marc. {\it
Produit complet des groupes de permutations et le probl\`eme
d'extension des groupes III.} {\sl Acta Sci.\ Math.\ Szeged} {\bf 14}
(1951) pp.~\hbox{69--82}. {MR:14,242d}}\par\filbreak\fi
\ifundefined{xkhukhro}\else
\item{\bf [\refer{khukhro}]}{Khukhro, Evgenii I. {\it Nilpotent Groups
and their Automorphisms.} {\sl de Gruyter Expositions in Mathematics}
{\bf 8}, New York 1993. {MR:94g:20046}}\par\filbreak\fi
\ifundefined{xkleimanbig}\else
\item{\bf [\refer{kleimanbig}]}{Kle\u{\i}man, Yu.~G. {\it On
identities in groups.\/} {\sl Trans.\ Moscow Math.\ Soc.\ } 1983,
Issue 2, pp.~\hbox{63--110}. {MR:84e:20040}}\par\filbreak\fi
\ifundefined{xthirtynine}\else
\item{\bf [\refer{thirtynine}]}{Kov\'acs, L.~G. {\it The thirty-nine
varieties.} {\sl Math.\ Scientist} {\bf 4} (1979)
pp.~\hbox{113--128.} {MR:81m:20037}}\par\filbreak\fi
\ifundefined{xlamssix}\else
\item{\bf [\refer{lamssix}]}{Lam, T{.}Y{.}, and Leep, David B. {\it
Combinatorial structure on the automorphism group of~$S_6$.\/} {\sl
Expo. Math.} {\bf 11} (1993) pp.~\hbox{289--308.}
{MR:94i:20006}}\par\filbreak\fi
\ifundefined{xlevione}\else
\item{\bf [\refer{levione}]}{Levi, F.~W. {\it Groups on which the
commutator relation 
satisfies certain algebraic conditions.\/} {\sl J.\ Indian Math.\ Soc.\ New
Series} {\bf 6}(1942), pp.~\hbox{87--97.} {MR:4,133i}}\par\filbreak\fi
\ifundefined{xgermanlevi}\else
\item{\bf [\refer{germanlevi}]}{Levi, F.~W. and van der Waerden,
B.~L. {\it \"Uber eine 
besondere Klasse von Gruppen.\/} {\sl Abhandl.\ Math.\ Sem.\ Univ.\ Hamburg}
{\bf 9}(1932), pp.~\hbox{154--158.}}\par\filbreak\fi
\ifundefined{xlichtman}\else
\item{\bf [\refer{lichtman}]}{Lichtman, A.~L. {\it Necessary and
sufficient conditions for the residual nilpotence of free products of
groups.\/} {\sl J. Pure and Applied Algebra} {\bf 12} no. 1 (1978),
pp.~\hbox{49--64.} {MR:58\#5938}}\par\filbreak\fi
\ifundefined{xmaxofan}\else
\item{\bf [\refer{maxofan}]}{Liebeck, Martin W.; Praeger, Cheryl E.;
and Saxl, Jan. {\it A classification of the maximal subgroups of the
finite alternating and symmetric groups.\/} {\sl J. of Algebra} {\bf
111}(1987), pp.~\hbox{365--383.} {MR:89b:20008}}\par\filbreak\fi
\ifundefined{xepisingroups}\else
\item{\bf [\refer{episingroups}]}{Linderholm, C.E. {\it A group
epimorphism is surjective.\/} {\sl Amer.\ Math.\ Monthly\ }77
pp.~\hbox{176--177.}}\par\filbreak\fi
\ifundefined{xmckay}\else
\item{\bf [\refer{mckay}]}{McKay, Susan. {\it Surjective epimorphisms
in classes
of groups.} {\sl Quart.\ J.\ Math.\ Oxford (2),\/} {\bf 20} (1969),
pp.~\hbox{87--90.} {MR:39\#1558}}\par\filbreak\fi
\ifundefined{xmaclane}\else
\item{\bf [\refer{maclane}]}{Mac Lane, Saunders. {\it Categories for
the Working Mathematician.} {\sl Graduate texts in mathematics 5},
Springer Verlag (1971). {MR:50\#7275}}\par\filbreak\fi
\ifundefined{xbilinear}\else
\item{\bf [\refer{bilinear}]}{Magidin, Arturo. {\it Bilinear maps and 
2-nilpotent groups.\/} August 1996, 7~pp.}\par\filbreak\fi
\ifundefined{xbilinearprelim}\else
\item{\bf [\refer{bilinearprelim}]}{Magidin, Arturo. {\it Bilinear maps
and central extensions of abelian groups.\/} In~preparation.}\par\filbreak\fi
\ifundefined{xprodvar}\else
\item{\bf [\refer{prodvar}]}{Magidin, Arturo. {\it Dominions in product
varieties of groups.\/} May 1997, 21~pp.}\par\filbreak\fi
\ifundefined{xprodvarprelim}\else
\item{\bf [\refer{prodvarprelim}]}{Magidin, Arturo. {\it Dominions in product
varieties of groups.\/} In preparation.}\par\filbreak\fi
\ifundefined{xmythesis}\else
\item{\bf [\refer{mythesis}]}{Magidin, Arturo. {\it Dominions in
Varieties of Groups.\/} Doctoral dissertation, University of
California at Berkeley, May 1998.}\par\filbreak\fi
\ifundefined{xnildoms}\else
\item{\bf [\refer{nildoms}]}{Magidin, Arturo {\it Dominions in varieties
of nilpotent groups.\/} December 1996, 27~pp.}\par\filbreak\fi
\ifundefined{xnildomsprelim}\else
\item{\bf [\refer{nildomsprelim}]}{Magidin, Arturo. {\it Dominions in
varieties of nilpotent groups.\/} In preparation.}\par\filbreak\fi
\ifundefined{xsimpleprelim}\else
\item{\bf [\refer{simpleprelim}]}{Magidin, Arturo. {\it Dominions in
varieties generated by simple groups.\/} In preparation.}\par\filbreak\fi
\ifundefined{xntwodoms}\else
\item{\bf [\refer{ntwodoms}]}{Magidin, Arturo. {\it Dominions in the variety of
2-nilpotent groups.\/} May 1996, 6~pp.}\par\filbreak\fi
\ifundefined{xdomsmetabprelim}\else
\item{\bf [\refer{domsmetabprelim}]}{Magidin, Arturo. {\it Dominions
in the variety of metabelian groups.\/}
In~preparation.}\par\filbreak\fi
\ifundefined{xfgnilgroups}\else
\item{\bf [\refer{fgnilgroups}]}{Magidin, Arturo. {\it Dominions of
finitely generated nilpotent groups.\/} October~1997,
10~pp.}\par\filbreak\fi
\ifundefined{xfgnilprelim}\else
\item{\bf [\refer{fgnilprelim}]}{Magidin, Arturo. {\it Dominions of
finitely generated nilpotent groups.\/} In preparation.}\par\filbreak\fi
\ifundefined{xwordsprelim}\else
\item{\bf [\refer{wordsprelim}]}{Magidin, Arturo. {\it
Words and dominions.\/} In~preparation.}\par\filbreak\fi
\ifundefined{xepis}\else
\item{\bf [\refer{epis}]}{Magidin, Arturo. {\it Non-surjective epimorphisms
in varieties of groups and other results.\/} February 1997,
13~pp.}\par\filbreak\fi
\ifundefined{xoddsandends}\else
\item{\bf [\refer{oddsandends}]}{Magidin, Arturo. {\it Some odds and
ends.\/} June 1996, 3~pp.}\par\filbreak\fi
\ifundefined{xpropdom}\else
\item{\bf [\refer{propdom}]}{Magidin, Arturo. {\it Some properties of
dominions in varieties of groups.\/} March 1997, 13~pp.}\par\filbreak\fi
\ifundefined{xzabsp}\else
\item{\bf [\refer{zabsp}]}{Magidin, Arturo. {\it $\Z$ is an absolutely
closed $2$-nil group.\/} Submitted.}\par\filbreak\fi
\ifundefined{xmagnus}\else
\item{\bf [\refer{magnus}]}{Magnus, Wilhelm; Karras, Abraham; and
Solitar, Donald. {\it Combinatorial Group Theory.\/} 2nd Edition; Dover
Publications, Inc.~1976. {MR:53\#10423}}\par\filbreak\fi
\ifundefined{xamalgtwo}\else
\item{\bf [\refer{amalgtwo}]}{Maier, Berthold J. {\it Amalgame
nilpotenter Gruppen
der Klasse zwei II.\/} {\sl Publ.\ Math.\ Debrecen} {\bf 33}(1986),
pp.~\hbox{43--52.} {MR:87k:20050}}\par\filbreak\fi
\ifundefined{xnilexpp}\else
\item{\bf [\refer{nilexpp}]}{Maier, Berthold J. {\it On nilpotent
groups of exponent $p$.\/} {\sl Journal of~Algebra} {\bf 127} (1989)
pp.~\hbox{279--289.} {MR:91b:20046}}\par\filbreak\fi
\ifundefined{xmaltsev}\else
\item{\bf [\refer{maltsev}]}{Maltsev, A.~I. {\it Generalized
nilpotent algebras and their associated groups.} (Russian) {\sl
Mat.\ Sbornik N.S.} {\bf 25(67)} (1949) pp.~\hbox{347--366.} ({\sl
Amer.\ Math.\ Soc.\ Translations Series 2} {\bf 69} 1968,
pp.~\hbox{1--21.}) {MR:11,323b}}\par\filbreak\fi
\ifundefined{xmaltsevtwo}\else
\item{\bf [\refer{maltsevtwo}]}{Maltsev, A.~I. {\it Homomorphisms onto
finite groups.} (Russian) {\sl Ivanov. gosudarst. ped. Inst., u\v
cenye zap., fiz-mat. Nuak} {\bf 18} (1958)
\hbox{pp. 49--60.}}\par\filbreak\fi
\ifundefined{xmorandual}\else
\item{\bf [\refer{morandual}]}{Moran, S. {\it Duals of a verbal
subgroup.\/} {\sl J.\ London Math.\ Soc.} {\bf 33} (1958)
pp.~\hbox{220--236.} {MR:20\#3909}}\par\filbreak\fi
\ifundefined{xhneumann}\else
\item{\bf [\refer{hneumann}]}{Neumann, Hanna. {\it Varieties of
Groups.\/} {\sl Ergebnisse der Mathematik und ihrer Grenz\-ge\-biete\/}
New series, Vol.~37, Springer Verlag~1967. {MR:35\#6734}}\par\filbreak\fi
\ifundefined{xneumannwreath}\else
\item{\bf [\refer{neumannwreath}]}{Neumann, Peter M. {\it On the
structure of standard wreath products of groups.\/} {\sl Math.\
Zeitschr.\ }{\bf 84} (1964) pp.~\hbox{343--373.} {MR:32\#5719}}\par\filbreak\fi
\ifundefined{xpneumann}\else
\item{\bf [\refer{pneumann}]}{Neumann, Peter M. {\it Splitting groups
and projectives
in varieties of groups.\/} {\sl Quart.\ J.\ Math.\ Oxford} (2), {\bf
18} (1967),
pp.~\hbox{325--332.} {MR:36\#3859}}\par\filbreak\fi
\ifundefined{xoates}\else
\item{\bf [\refer{oates}]}{Oates, Sheila. {\it Identical Relations in
Groups.\/} {\sl J.\ London Math.\ Soc.} {\bf 38} (1963),
pp.~\hbox{71--78.} {MR:26\#5043}}\par\filbreak\fi
\ifundefined{xolsanskii}\else
\item{\bf [\refer{olsanskii}]}{Ol'\v{s}anski\v{\i}, A. Ju. {\it On the
problem of a finite basis of identities in groups.\/} {\sl
Izv.\ Akad.\ Nauk.\ SSSR} {\bf 4} (1970) no. 2
pp.~\hbox{381--389.}}\par\filbreak\fi
\ifundefined{xremak}\else
\item{\bf [\refer{remak}]}{Remak, R. {\it \"Uber minimale invariante
Untergruppen in der Theorie der end\-lichen Gruppen.\/} {\sl
J.\ reine.\ angew.\ Math.} {\bf 162} (1930),
pp.~\hbox{1--16.}}\par\filbreak\fi
\ifundefined{xclassifthree}\else
\item{\bf [\refer{classifthree}]}{Remeslennikov, V. N. {\it Two
remarks on 3-step nilpotent groups} (Russian) {\sl Algebra i Logika
Sem.} (1965) no.~2 pp.~\hbox{59--65.} {MR:31\#4838}}\par\filbreak\fi
\ifundefined{xrotman}\else
\item{\bf [\refer{rotman}]}{Rotman, J.J. {\it Introduction to the Theory of
Groups}, 4th edition. {\sl Graduate texts in mathematics 119},
Springer Verlag,~1994. {MR:95m:20001}}\par\filbreak\fi
\ifundefined{xsaracino}\else
\item{\bf [\refer{saracino}]}{Saracino, D. {\it Amalgamation bases for
nil-$2$ groups.\/} {\sl Alg.\ Universalis\/} {\bf 16} (1983),
pp.~\hbox{47--62.} {MR:84i:20035}}\par\filbreak\fi
\ifundefined{xscott}\else
\item{\bf [\refer{scott}]}{Scott, W.R. {\it Group Theory.} Prentice
Hall,~1964. {MR:29\#4785}}\par\filbreak\fi
\ifundefined{xsmelkin}\else
\item{\bf [\refer{smelkin}]}{\v{S}mel'kin, A.L. {\it Wreath products and
varieties of groups} [Russian] {\sl Dokl.\ Akad.\ Nauk S.S.S.R.\/} {\bf
157} (1964), pp.~\hbox{1063--1065} Transl.: {\sl Soviet Math.\ Dokl.\ } {\bf
5} (1964), pp.~\hbox{1099--1011}. {MR:33\#1352}}\par\filbreak\fi
\ifundefined{xstruikone}\else
\item{\bf [\refer{struikone}]}{Struik, Ruth Rebekka. {\it On nilpotent
products of cyclic groups.\/} {\sl Canadian Journal of
Mathematics\/} {\bf 12} (1960)
pp.~\hbox{447--462}. {MR:22\#11028}}\par\filbreak\fi
\ifundefined{xstruiktwo}\else
\item{\bf [\refer{struiktwo}]}{Struik, Ruth Rebekka. {\it On nilpotent
products of cyclic groups II.\/} {\sl Canadian Journal of
Mathematics\/} {\bf 13} (1961) pp.~\hbox{557--568.}
{MR:26\#2486}}\par\filbreak\fi
\ifundefined{xvlee}\else
\item{\bf [\refer{vlee}]}{Vaughan-Lee, M{.} R{.} {\it Uncountably many
varieties of groups.\/} {\sl Bull.\ London Math.\ Soc.} {\bf 2} (1970)
pp.~\hbox{280--286.} {MR:43\#2054}}\par\filbreak\fi
\ifundefined{xweibel}\else
\item{\bf [\refer{weibel}]}{Weibel, Charles. {\it Introduction to
Homological Algebra.\/} Cambridge University
Press~1994. {MR:95f:18001}}\par\filbreak\fi 
\ifundefined{xweigelone}\else
\item{\bf [\refer{weigelone}]}{Weigel, T.S. {\it Residual properties
of free groups.\/} {\sl J.\ of Algebra} {\bf 160} (1993)
pp.~\hbox{14--41.} {MR:94f:20058a}}\par\filbreak\fi
\ifundefined{xweigeltwo}\else
\item{\bf [\refer{weigeltwo}]}{Weigel, T.S. {\it Residual properties
of free groups II.\/} {\sl Comm.\ in Algebra} {\bf 20}(5) (1992)
pp.~\hbox{1395--1425.} {MR:94f:20058b}}\par\filbreak\fi
\ifundefined{xweigelthree}\else 
\item{\bf [\refer{weigelthree}]}{Weigel, T.S. {\it Residual Properties
of free groups III.\/} {\sl Israel J.\ Math.\ } {\bf 77} (1992)
pp.~\hbox{65--81.} {MR:94f:20058c}}\par\filbreak\fi
\ifundefined{xzstwo}\else
\item{\bf [\refer{zstwo}]}{Zariski, Oscar and Samuel, Pierre. {\it
Commutative Algebra}, Volume
II. Springer-Verlag~1976. {MR:52\#10706}}\par\filbreak\fi
\ifnum0<\citations\nonfrenchspacing\fi

\bigskip
{\it
\obeylines
\noindent Arturo Magidin
\noindent Cub\'iculo 112
\noindent Instituto de Matem\'aticas
\noindent Universidad Nacional Aut\'onoma de M\'exico
\noindent 04510 Mexico City, MEXICO
\noindent e-mail: magidin@matem.unam.mx
}
 
\vfill\eject
\immediate\closeout\aux
\end